\documentclass[a4paper,11pt]{article}
\usepackage{amsmath, amsfonts, amssymb, amsthm, graphicx}
\usepackage{amsmath,amssymb}
\usepackage{abstract}
\usepackage{amssymb}
\usepackage{amsmath}
\usepackage{graphics}
\usepackage{relsize}
\usepackage{latexsym}
\usepackage{amsfonts}
\usepackage{color}
\usepackage{epstopdf}
\usepackage{ctable}
\usepackage{graphicx}
\usepackage{array,multirow,makecell}
\usepackage{verbatim}
\newtheorem{theo}{Theorem}[section]

\newtheorem{prop}[theo]{Proposition}
\newtheorem{pf}[theo]{Proof} 
\newtheorem{uppg}[theo]{Exercise}
\newtheorem{remark}[theo]{Remark}
\newcommand{\be}{\begin{eqnarray*}}
\newcommand{\ee}{\end{eqnarray*}}
\newcommand{\ben}{\begin{eqnarray}}
\newcommand{\een}{\end{eqnarray}}

\def\subsecn (#1) {\medskip\ \ \ {\it #1}\medskip}

\newcommand{\lp}[1]{\left(\begin{array}{#1}}
\newcommand{\rp}{\end{array}\right)}

\newcommand{\leftd}[1]{\left\{\begin{array}{#1}}
\newcommand{\rightd}{\end{array}\right.}
\setcellgapes{1pt}
\makegapedcells
\newcolumntype{R}[1]{>{\raggedleft\arraybackslash }b{#1}}
\newcolumntype{L}[1]{>{\raggedright\arraybackslash }b{#1}}
\newcolumntype{C}[1]{>{\centering\arraybackslash }b{#1}}

\def\A {\mathbf{A}}

\def\I {\mathbf{I}}

\def\R {\mathbf{R}}

\def\l {\boldsymbol{l}}

\def\s {\boldsymbol{s}}

\def\w {\boldsymbol{w}}
\def\x {\boldsymbol{x}}
\def\y {\boldsymbol{y}}

\def\Bc {\mathcal{B}}

\def\Nc {\mathcal{N}}

\def\Rb {\mathbb{R}}
\def\Pb {\mathbb{P}}

\begin{document}

\begin{center} 
\textbf{\large MCMC convergence diagnosis using geometry of   
Bayesian LASSO}\\
\end{center}
\begin{center} 
A. Dermoune, D.Ounaissi, N.Rahmania
\end{center}





\begin{abstract}
Using posterior distribution of Bayesian LASSO  
we construct a semi-norm on the parameter space.  
We show that the partition function depends on the ratio
of the $l^1$ and $l^2$ norms and present three regimes.  
We derive the concentration of Bayesian LASSO, and 
present MCMC convergence diagnosis. 

\end{abstract}

\textbf{keyword:}
LASSO, Bayes, MCMC, log-concave, geometry, incomplete Gamma function

\section{Introduction} 
Let $p \geq n$ be two positive integers, $\y \in \Rb^n$ and 
$\A$ be an $n\times p$ matrix with real numbers entries. Bayesian LASSO 
\ben 
c(\x)=\frac{1}{Z}\exp\Big(-\frac{\|\A\x-\y\|_2^2}{2}-\|\x\|_1\Big)\label{c}
\een   
is a typically posterior distribution used in the linear regression 
\be 
\y=\A\x+\w. 
\ee 
Here 
\ben 
Z=\int_{\Rb^p}\exp\Big(-\frac{\|\A\x-\y\|_2^2}{2}-\|\x\|_1\Big)d\x\label{Zx}
\een

is the partition function, $\|\cdot\|_2$ and $\|\cdot\|_1$ are respectively 
the Euclidean and the $l_1$ norms. 
The vector $\y\in\Rb^n$ are the observations, $\x\in \Rb^p$ is the unknown signal to recover, $\w\in\Rb^n$ is the standard Gaussian noise, and $\A$ is a known matrix which maps the signal domain $\Rb^p$ into the observation domain $\Rb^n$. If we suppose that $\x$ is drawn from 
Laplace distribution i.e. the distribution proportional to 
\ben 
\exp(-\|\x\|_1),
\label{laplace}  
\een
then the posterior of $\x$ known $y$ is drawn from the distribution $c$ (\ref{c}). 
The mode  
\ben 
\arg\min\Big\{\frac{\|\A\x-\y\|_2^2}{2}+\|\x\|_1:\quad \x\in\Rb^p\Big\}
\een
of $c$ was first introduced in \cite{Tibshirani1996} and called LASSO. It is also called Basis Pursuit De-Noising method \cite{Donoho}. In our work we select the term LASSO and keep it for the rest of the article.  

In general LASSO is not a singleton, i.e. the mode of the distribution $c$ 
is not unique. In this case LASSO is a set and we will denote by lasso any element of this set. 
A large number of theoretical results has been provided 
for LASSO. See \cite{Daubechies2004}, \cite{DDN}, \cite{Fort}, \cite{Pereyra} and the references herein. 
The most popular algorithms to find LASSO are LARS algorithm \cite{Efron}, ISTA and FISTA algorithms see e.g. \cite{Beck} and the review article \cite{Parikh}.

The aim of this work is to study geometry of bayesian LASSO 
and to derive MCMC convergence diagnosis. 
\section{Polar integration}   
Using polar coordinates, the partition function (\ref{Zx})  
\ben 
Z=\int_{S}J_p(\theta)d\theta,
\label{zpolaire} 
\een 
where $d\theta$ denotes the surface measure on the unit sphere $S$, and 
\ben 
J_p(\theta)=\int_0^{+\infty}\exp(-g(r,\theta))r^{p-1}dr.  
\label{Jp}
\een 
Here 
\ben 
g(r,\theta)=\frac{1}{2}(r^2\|\A\theta\|_2^2+2r\|\A\theta\|_2\beta
+\|\y\|_2^2),  
\label{g} 
\een
where 
\ben 
\beta:=\frac{\|\theta\|_{1}}{\|\A\theta\|_2}-\|\y\|_2s,  
\label{beta}
\een 
and $s$ denotes the cosine of the angle $(\A\theta, \y)$ i.e. $\cos((\A\theta, \y))$ . Using known estimate $\|\theta\|_2\leq \|\theta\|_1$, we observe that $\beta$ is bounded below by 
\ben 
\frac{1}{\|\A\|}-\|\y\|_2, 
\label{betalowerbound} 
\een 
and 
$\beta\to +\infty$ as $\A\theta\to 0$. Here 
$\|\A\|$ is the square root of the largest eigenvalue of $\A^*\A$.   
Observe that 
\be 
c(\x)d\x=\frac{1}{Z}\exp(-g(r,\theta))r^{p-1}drd\theta.  
\ee
Hence, we can sample from Bayesian LASSO $c$ (\ref{c}) as following. We draw  
uniformly $\theta:=\frac{\x}{\|\x\|_2}$ 
from the unit sphere, and then draw the norm $\|\x\|_2$ following the distribution 
\ben 
\mu_{\theta}(r):=\frac{1}{J_p(\theta)}\exp(-\varphi(r,\theta)), 
\label{mutheta} 
\een 
where   
\ben 
\varphi(r,\theta):=g(r,\theta)-(p-1)\ln(r),\quad r >0.  
\label{varphi} 
\een
Moreover, observe that the modes $\{\x_{lasso}=r_{lasso}\theta_{lasso}:\quad 
g(r_{lasso},\theta_{lasso})=
\min_{r\geq 0, \theta\in S}g(r,\theta)\}$ and  
$\{(r^*,\theta^*):\quad 
\varphi(r,\theta)=\min_{r\geq 0, \theta\in S}\varphi(r,\theta)\}$
respectively of the distributions $c(\x)d\x$ and $\frac{1}{Z}\exp(-g(r,\theta))r^{p-1}drd\theta$ 
are different. We will show that $(r^*,\theta^*)$ contains more information 
than $\x_{lasso}$. 
 
\section{Geometric interpretation of the partition function} 
The volume (Lebesgue measure) of the set $K(\A,\y):=\{\x\in\Rb^p: J_p(\x)\geq 1\}$ is 
$vol(K(\A,\y))=\frac{1}{p} Z$.  
Observe that $J_p^{-\frac{1}{p}}$ is a norm on the null-space $N(\A)$ of $\A$. 
A general result \cite{Ball} tells us that if $f$ is even, log-concave and 
integrable on an Euclidean space $E$, then 
\be 
\x\in E\to (\int_0^{+\infty}f(r\x)r^{p-1}dr)^{-\frac{1}{p}}
\ee
is a norm on $E$. It follows that in the case $E=N(\A)$ or $E=\{\x\in\Rb^p: 
\langle \A\x,\y\rangle=0\}$, the map 
\be 
\x\in E\to J_p^{-\frac{1}{p}}(\x)
\ee 
is a norm on $E$. The map 
\be 
\x\in\Rb^p\to J_p^{-\frac{1}{p}}(\x):=\|\x\|_{LASSO}
\ee 
has nearly all the properties of a norm. Only the evenness is missing.
The set $K(\A,\y)=\{\x\in\Rb^p: \quad \|\x\|_{LASSO}\leq 1\}$ 
is convex, compact and contains the origin. See \cite{Klartag} 
for more details. 

\subsection{Necessary and sufficient condition to have $LASS0=\{0\}$} 
If $\beta \geq 0$, then $r\in [0,+\infty)\to g(r,\theta)$ 
is increasing, its minimizer is equal to $r=0$, and  
its smallest value is $\frac{\|\y\|_2^2}{2}$. If $\beta <0$, then 
its minimizer is equal to $r=-\frac{\beta}{\|\A\theta\|_2}$, and  
its smallest value is less than $\frac{\|\y\|_2^2}{2}$. 
If the set $\{\beta < 0\}$ is empty, then 
$LASSO=\{0\}$, if not 
\ben 
LASSO=\{\l=-\frac{\beta_l}{\|\A\theta_l\|_2}\theta_l:\quad \beta_l\leq 0,\,\mbox{s.t.}\, 
\beta_l^2=\sup_{\beta \leq 0}\beta^2\}.
\label{polarlasso} 
\een 
As an illustration we consider the case $n=4$, $p=7$ 
and the entries of the matrix $\A \sim \Bc(\pm\frac{1}{\sqrt{n}})$
are a realization of i.i.d. Bernoulli random variables with the values $\pm\frac{1}{\sqrt{n}}$. We draw uniformly $N=10^5$ 
vectors from the sphere $S$ and estimate LASSO using Formula (\ref{polarlasso}). 
Table 1 gives the value of LASSO using respectively  
FISTA algorithm and Formula (\ref{polarlasso}).   
\renewcommand{\arraystretch}{0.9} 
\setlength{\tabcolsep}{0.07cm} 
 \begin{table}[!ht]
\begin{center}
\begin{tabular}{|c|c|c|c|c|c|c|c|}
\hline  &$x_1$ &
 $ x_2$ &
   $x_3$&
   $x_4$&
   $x_5$&
   $x_6$&
   $x_7$\\
\hline $LASSO_{FISTA}$ &1.7744&
  0.6019 &
   -0.3283 &
   0 &
    0&
   -1.0050&
  0 \\
\hline  
 $LASSO_{POLAR}$& 0.9992 &
 0.3890&
 -1.3980&
  0.0769&
   -0.0070&
  -0.8699&
   -0.0379\\
\hline 
\end{tabular}
\caption{ $N=10^5$, $p=7$, $n=4$.}
\end{center}
\end{table}

Observe that necessary and sufficient condition for $\beta\geq 0$ for all $\theta$ is 
\be 
\|\y\|_2\leq \inf\{\frac{\|\theta\|_1}{\|\A\theta\|_2|s|}:\quad \theta\in S\}.
\ee  
Using known estimate $\|\theta\|_2\leq \|\theta\|_1$, we obtain
\be 
\|\y\|_2\leq \frac{1}{\|\A\|}
\ee 
as a sufficient condition for $\beta \geq 0$ for all $\theta$. 

\section{Closed form of the partition function} 

We introduce for $a\in \Rb$ 
and for a couple $p, r\geq 1$ of integers, the notations 
\ben 
(a)_r&=&(a-1)\ldots (a-r),\\
c(p,r)&=&\sum_{k=0}^{p-1}\binom{p-1}{k}(-1)^{p-1-k}
(\frac{k+1}{2})_r.
\label{cpr}
\een 
Now, we can announce the following result.  

\begin{prop} \label{closedfromlassozero} 1) If $\beta=+\infty$, then 
\be 
\|\theta\|_1^pJ_p(\theta)=(p-1)!\exp(-\frac{\|\y\|_2^2}{2}).
\ee 

2) If $\beta \geq 0$, then 
\be
\|\theta\|_1^pJ_{p}(\theta):=\Phi(\beta),
\ee 
where 
\ben 
\Phi(\beta)=&\exp(-\frac{\|\y\|_2^2}{2})(\beta+s\|\y\|_2)^p
\sum_{k=0}^{p-1}\binom{p-1}{k}(-\beta)^{p-1-k}\nonumber\\
&2^{\frac{k-1}{2}}\exp(\frac{\beta^2}{2})\Gamma(\frac{k+1}{2},\frac{\beta^2}{2})\label{phibeta}, 
\een 
Here $\Gamma(a,x)=\int_x^{+\infty}\exp(-t)t^{a-1}dt,\quad a>0, x\geq 0$, is the upper incomplete Gamma function.\\

3) If $\beta < 0$, then 
\ben 
\|\theta\|_1^pJ_p(\theta)=&\exp(-\frac{\|\y\|_2^2}{2})(\beta+s\|\y\|_2)^p
\sum_{k=0}^{p-1}\binom{p-1}{k}(-\beta)^{p-1-k}\nonumber\\
&2^{\frac{k-1}{2}}\exp(\frac{\beta^2}{2})2^{\frac{k-1}{2}}(\Gamma(\frac{k+1}{2})+(-1)^k\gamma(\frac{k+1}{2},\frac{\beta^2}{2})).  
\label{betanegatifformula} 
\een  
Here $\gamma(a,x)=\int_0^x\exp(-t)t^{a-1}dt$ is the lower incomplete Gamma function.\\ 

4) If $\beta=0$, then 
\be 
\|\theta\|_1^pJ_p(\theta)=\exp(-\frac{\|\y\|_2^2}{2})2^{\frac{p-2}{2}}\|\y\|_2^ps^p\Gamma(\frac{p}{2},0).
\ee 

5) If $\beta >0$ then for $M\geq p+1$, 
\be
\|\theta\|_1^pJ_p(\theta):=\Phi(\beta,M)+R(\beta,M),
\ee 
where 
\ben  
\Phi(\beta,M)=(p-1)!\exp(-\frac{\|\y\|_2^2}{2})+\sum_{r=p}^{M-1}
2^{p-1}\exp(-\frac{\|\y\|_2^2}{2})\Big(1+\frac{\|\y\|_2s}{\beta}\Big)^pc(p,r) 
\Big(\frac{\beta^2}{2}\Big)^{p-1-r}\label{phibetam} 
\een
and the remainder term 
\ben 
 |R(\beta,M)|\leq \Big(1+\frac{\|\y\|_2s}{\beta}\Big)^p\exp(-\frac{\|\y\|_2^2}{2})
2^{p-1}\frac{|c(p,M)|}{\Big(\frac{\beta^2}{2}\Big)^{M-(p-1)}}\label{RT} 
\een 
\end{prop}

\begin{pf} Only the first part of assertions 2) and 5) needs the proof. 
Let us prove the first part of 2). From the equality 
\be
J_{p}(\theta)=\exp(-\frac{\|\y\|_2^2}{2}+\frac{\beta^2}{2})\int_0^{+\infty}\exp(-\frac{(\|\A\theta\|_2r+\beta)^2}{2})r^{p-1}dr, 
\ee 
and the change of the variable 
\be 
\tau=\|\A\theta\|_2r+\beta,
\ee 
we obtain 
\be 
\int_0^{+\infty}\exp\Big(-\frac{(\|\A\theta\|_2r+\beta)^2}{2}\Big)r^{p-1}dr&=&\frac{1}{\|\A\theta\|_2^p}\int_{\beta}^{+\infty}\exp(-\frac{\tau^2}{2})(\tau-\beta)^{p-1}d\tau \\
&=&\frac{1}{\|\A\theta\|_2^p}\sum_{k=0}^{p-1}\binom{p-1}{k}(-\beta)^{p-1-k}
\int_{\beta}^{+\infty}\exp(-\frac{\tau^2}{2})\tau^kd\tau. 
\ee 
As $\beta >0$, then the change of variable 
\be 
\omega=\frac{\tau^2}{2}, 
\ee 
implies 
\be 
\int_{\beta}^{+\infty}\exp(-\frac{\tau^2}{2})\tau^kd\tau=2^{\frac{k-1}{2}}\Gamma(\frac{k+1}{2},\frac{\beta^2}{2}). 
\ee  
The equality $\beta=\frac{\|\theta\|_1}{\|\A\theta\|_2}-\|\y\|_2s$ achieves the proof of 2).  

Now we prove Assertion 5). We extend the incomplete Gamma function 
as following 
\be 
\Gamma(a,x)=\int_x^{+\infty}\exp(-t)t^{a-1}dt,\quad x >0, a\in\Rb,
\ee 
and we use known estimate see \cite{Cops} page 14  
\ben 
\Gamma(a,x)=\exp(-x)x^{a-1}+
\sum_{r=1}^{M-1}(a)_r\exp(-x)x^{a-1-r}+R(a,x,M),
\label{gammaexpansion}
\een 
where $M >a-1$, $x >0$, $a\in\Rb$, and the remainder term  
\be 
|R(a,x,M)|\leq (a)_{M-1}\exp(-x)x^{a-1-M}.
\ee 
If $\beta >0$, then from $\beta=\frac{\|\theta\|_1}{\|\A\theta\|_2}-\|\y\|_2s$, we have  
\be 
J_p(\theta)=&\frac{2^{p-1}}{\|\theta\|_1^p}\exp(-\frac{\|\y\|_2^2}{2})\Big(1+\frac{\|\y\|_2s}{\beta}\Big)^p
\Big(\frac{\beta^2}{2}\Big)^{p-1}
\sum_{k=0}^{p-1}\binom{p-1}{k}(-1)^{p-1-k}\\
&\exp(\frac{\beta^2}{2})\Big(\frac{\beta^2}{2}\Big)^{-\frac{k-1}{2}}
\Gamma(\frac{k+1}{2},\frac{\beta^2}{2}).
\ee 
Using the expansion (\ref{gammaexpansion}), and the fact that 
$J_p(\theta)\to  \frac{(p-1)!}{\|\theta\|_1^p}\exp(-\frac{\|\y\|_2^2}{2})$ 
as $\beta\to +\infty$, we obtain for $r< p-1$
\be 
c(p,r)=0,\quad 1\leq r < p-1,\quad c(p,p-1)=\frac{(p-1)!}{2^{p-1}}. 
\ee
It follows the following expansion:  
\be 
\Big(1+\frac{\|\y\|_2s}{\beta}\Big)^{-p}\|\theta\|_1^p\exp(\frac{\|\y\|_2^2}{2})J_p(\theta)=
(p-1)!+2^{p-1}\sum_{r=p}^{M-1}\frac{c(p,r)}{\Big(\frac{\beta^2}{2}\Big)^{r-(p-1)}}+\tilde{R}(\beta,M),
\ee 
where the remainder term 
\be 
|\tilde{R}(\beta,M)|\leq 2^{p-1}\frac{|c(p,M)|}{\Big(\frac{\beta^2}{2}\Big)^{M-(p-1)}},
\ee 
which achieves the proof. 
\end{pf}

\subsection{Numerical calculations} 
As an illustration we consider $n=4$, $p=7$, $\A \sim \Bc(\pm\frac{1}{\sqrt{n}})$, and $\y=0$. 
The choice $M=17$ corresponds to the relative error 
$\frac{|R(\beta,M)|}{|\Phi(\beta,M)|}\leq 10^{-4}$ for $\beta\geq 7.5$. 
\begin{figure}[!ht]\centering
\includegraphics[width=14cm]{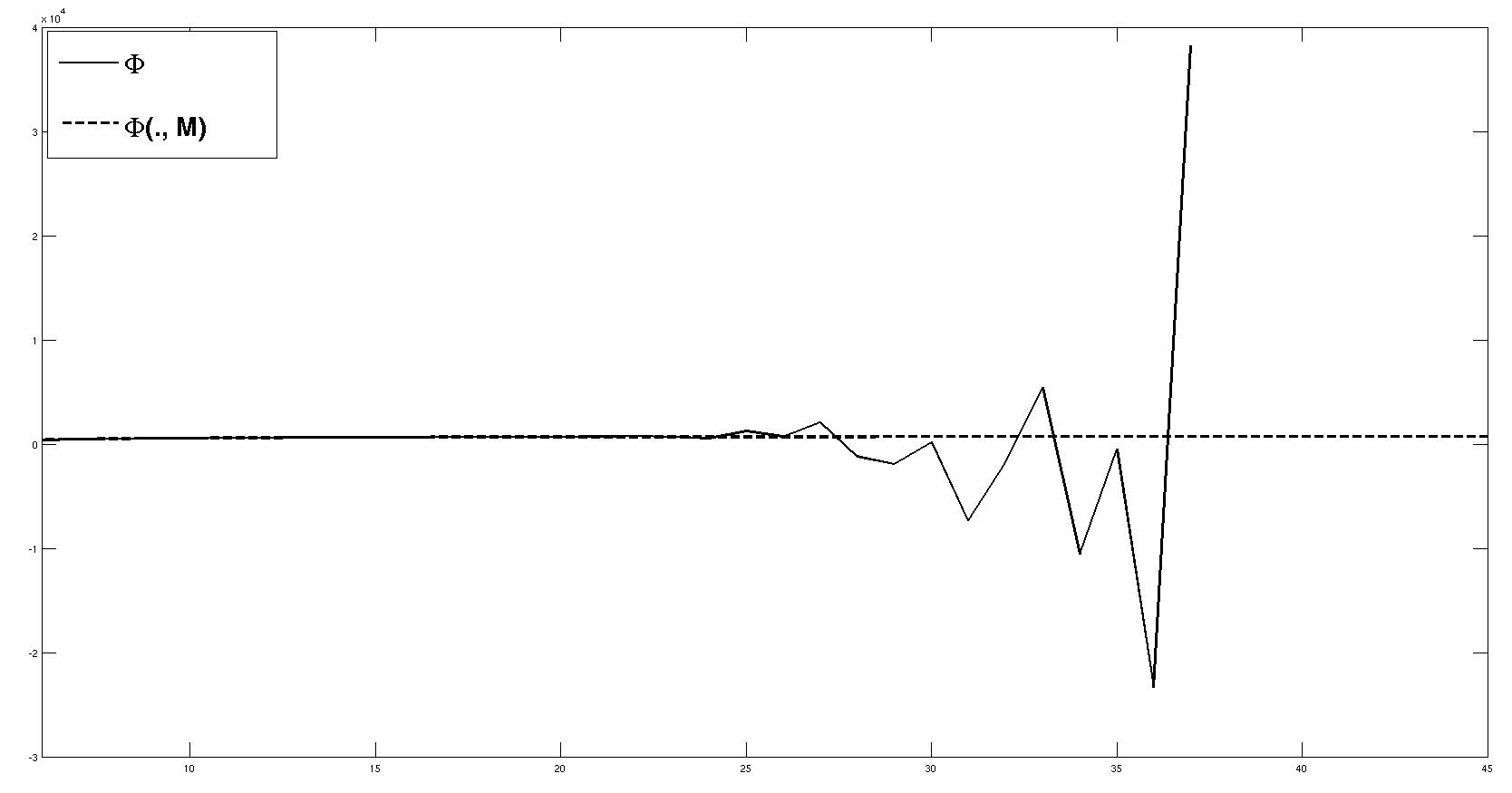} 
\caption{Curves of $\Phi(\beta)$ and $\Phi(\beta,17)$, $\beta=[6, 45]$, $p=7$, $n=4$.}\label{Fig1}
\end{figure}

\begin{figure}[!ht]\centering
\includegraphics[width=14.5cm]{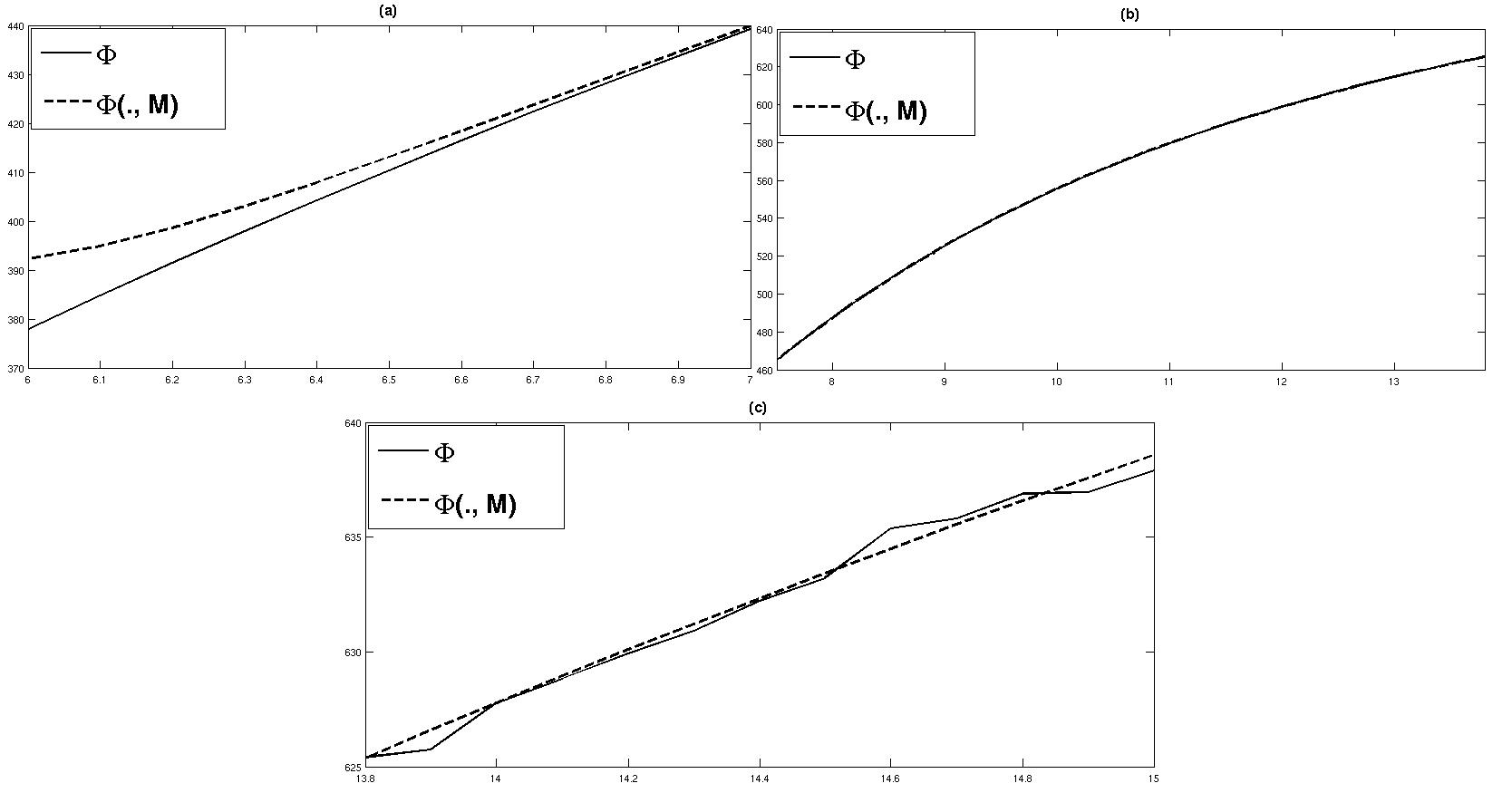} 
\caption{Curves of $\Phi(\beta)$ and $\Phi(\beta,17)$, $p=7$ and $n=4$.} \label{Fig2}
\end{figure}

Numerical calculations show that the function $\Phi(\beta)$ (\ref{phibeta}) 
explodes for $\beta > 13.8$. To compass 
these explosions we use the expansion (\ref{gammaexpansion}), and 
then we use the function $\Phi(\cdot,M)$ (\ref{phibetam}). 
In Fig.(\ref{Fig1}) and Fig.(\ref{Fig2}) dashed and black curves represent respectively the function $\Phi(\cdot,M)$ and 
$\Phi$. In Fig. (\ref{Fig1}) we plot $\Phi(\cdot,M)$ and 
$\Phi$  for $\beta\in (6, 45)$. By zooming on $\beta\in (1.09, 7.5)$, $\beta\in (7.5, 13.8)$, and  $\beta\in (13.8, 15)$ we show that the behavior of $\Phi$ becomes abnormal 
from $\beta \approx \beta_{\Phi}=13.8$ and obtain  
Fig.(\ref{Fig2}). 
\subsection{The case LASS0=\{ 0 \}: 
Partition function estimate and concentration inequality} 

The following is a consequence of \cite{Klartag} Lemma 2.1. 

\begin{figure}[!ht]\centering
\includegraphics[width=12cm]{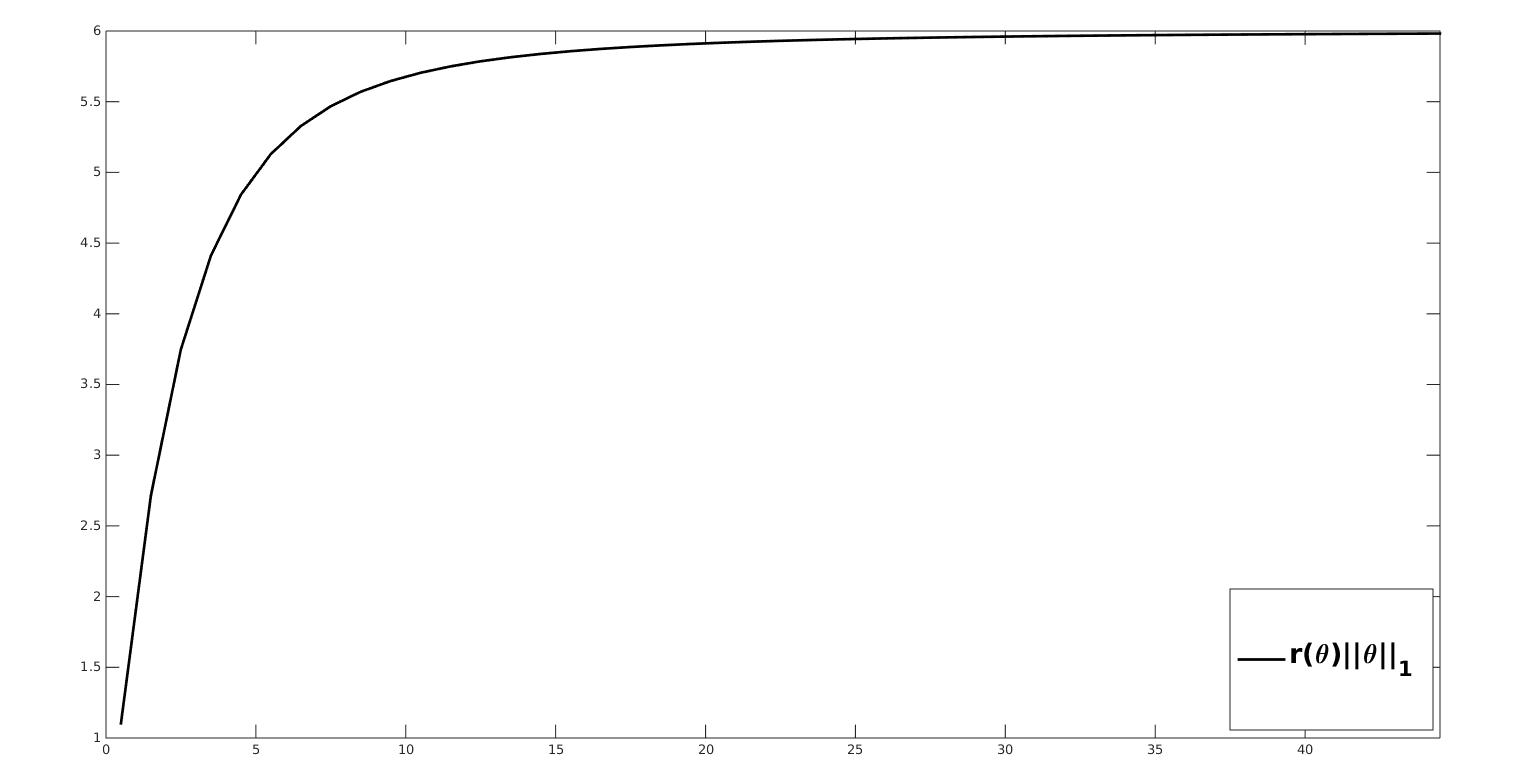} 
\caption{Curve of $r(\theta)\|\theta\|_1$, $\beta_{min}:= \frac{1}{\|A\|}= 0.4987$, $\beta_{max}=44.5$, $\min(r(\theta)\|\theta\|_1)= 1.1035$, $p=7$ and $n=4$.} \label{Fig3}
\end{figure}

\begin{prop} \label{lassozer} 1) The function  
$r\in (0, +\infty)\to \varphi(r,\theta)$
is convex and its unique critical point 
\ben 
r(\theta)=\frac{-\beta+\sqrt{\beta^2+4(p-1)}}{2\|\A\theta\|_2}
\label{rtheta} 
\een 
is the mode of (\ref{mutheta}).  

2) By denoting $M(\theta)=\exp\Big(-\varphi(r(\theta),\theta)\Big)$, we obtain 
 \be 
\frac{M(\theta)r(\theta)}{p}\leq  J_p(\theta)\leq \frac{M(\theta)r(\theta)(p-1)!\exp(p-1)}{(p-1)^{p}}. 
\ee

3) We have for $q>0$, 
\be 
\int_{qr(\theta)}\exp(-\varphi(r,\theta))dr\leq 
\frac{p\Gamma(p,(p-1)q)\exp(p-1)}{(p-1)^p}\int_0^{+\infty}\exp(-\varphi(r, \theta))dr,
\ee
and 
\be
Z_{min}:=\frac{|S|\inf_{\theta\in S}M(\theta)r(\theta)}{p}&\leq Z\leq \frac{|S|(p-1)!\exp(p-1)}{(p-1)^{p}}\sup_{\theta\in S}M(\theta)r(\theta):=Z_{max},
\ee 
where $|S|$ denotes the surface of the unit sphere $S$. 

4) If $\x$ is drawn from Bayesian LASSO distribution $c$ (\ref{c}), then for $q>0$, 
$\|\x\|_2\leq qr(\theta)$ with the probability at least equal to 
$P(q, p):= 1-\frac{p\Gamma(p,(p-1)q)\exp(p-1)}{(p-1)^p}$. In particular for $q=5$ we have 
$\frac{p\Gamma(p,(p-1)q)\exp(p-1)}{(p-1)^p}\leq \exp(-2(p-1))$.

\end{prop}
\begin{remark}
If $\y=0$, then LASS0=0 and the mode (\ref{rtheta}) becomes  
\be 
r(\theta)\|\theta\|_1=\frac{\beta(-\beta+\sqrt{\beta^2+4(p-1)}}{2}.\\
\ee   
\end{remark}
In Fig.(\ref{Fig3}) we plot $\beta\in [0.4987, 44.5]\to \frac{\beta(-\beta+\sqrt{\beta^2+4(p-1)}}{2}$.\\

The mode of the distribution of $\frac{1}{Z}\exp(-\varphi(r,\theta))drd\theta$ is equal to 
\be 
\arg\min_{r>0, \theta\in S}\varphi(r(\theta),\theta)=(r(\theta^*),\theta^*).
\ee 
As an illustration we consider $p=7$, $n=4$, $\A \sim \Bc(\pm \frac{1}{\sqrt{n}})$. 
We draw uniformly $N=10^5$ sample $\theta_i\in S$ from the unit sphere $S$. For 
each $i$, we calculate $\varphi(r(\theta_i),\theta_i)$, and we derive $\theta^*$. Notice that $\beta^*= \frac{\|\theta^*\|_1}{\|\A\theta^*\|_2}= 14.0122\approx \beta_{\Phi}$ is nearly equal to 
the beginning of abnormality of $\Phi$.

Using Formula (\ref{zpolaire}) and Monte Carlo method, we obtain  
$Z\approx 2.2142$, $Z_{min}\approx 0.0058$ and $Z_{max}\approx 120.3654$.  
If we draw $N=10^5$ vectors using Laplace distribution (\ref{laplace})
and calculate the value of $Z$ using Formula (\ref{Zx}) and Monte Carlo method, then 
we obtain $Z\approx 0.0036 < Z_{min}$. 
Hence Monte Carlo method using Formula (\ref{zpolaire}) wins against Monte Carlo 
method using Formula (\ref{Zx}). 
  

\section{The case $0\notin LASSO$}
If $0\notin LASSO$, then the assertions of Proposition (\ref{lassozer})
are no longer valid. But we are going to show that 
these assertions becomes valid if we work around LASSO. 
We consider for $\l\in LASSO$,
\ben 
h(\x)&=&-\frac{\|\A(\x+\l)-\y\|_2^2}{2}-\|\x+\l\|_1,\nonumber\\
\bar{h}(\x)&:=&h(\x)-h(0),\nonumber\\
f(\x)&=&\exp\Big(\bar{h}(\x)\Big)\label{f}. 
\een 
Contrary to the map $\x\to c(\x)$, the map 
$\x\to f(\x)$ attains its supremum at the origin. 
Observe that 
\be 
c(\x)&=&\frac{f(\x-\l)}{\int_{\Rb^p}f(\x)d\x},\\
Z&=&\exp(h(0))Z_f:=\exp(h(0))\int_{\Rb^p}f(\x)d\x.
\ee 
If $\x$ is drawn from $c$, then $\x-\l$ is drawn from $\frac{f}{Z_f}$. 
Moreover 
\ben 
Z_f=\int_{\Rb^p}f(\x)d\x=\int_SJ_{p}(\theta,\l)d\theta,
\label{zf}
\een 
where  
\ben 
J_p(\theta,\l):=\int_0^{+\infty}f(r\theta)r^{p-1}dr.  
\label{Jp}
\een 
The map $\x\in\Rb^p\to J_p^{-\frac{1}{p}}(\x,\l):=\|\x\|_{\A,\y,\l}$ 
is nearly a norm (only the eveness is missing). The set 
\ben 
K(\A,\y,\l)=\{\x\in\Rb^p: \|\x\|_{\A,\y,\l}\leq 1\}
\label{lassounitball}
\een 
is convex, compact and contains the origin. 
The volume 
\be 
Vol(K(\A,\y,\l))=\frac{Z_f}{p}. 
\ee 
If $\x$ is drawn from $c$, then $\x-\l$ is drawn from $\frac{f}{Z_f}$, or equivalently
if $\x$ is drawn from $\frac{f}{Z_f}$, then $\x+\l$ is drawn from $c$. 
To draw $\x$ from $\frac{f}{Z_f}$, we draw $\theta=\frac{\x}{\|\x\|_2}$ 
uniformely on $S$, and then we draw $\|\x\|_2$ from 
\ben 
\mu_{\theta,\l}(r):=\frac{f(r\theta)}{J_p(\theta,\l)}.
\label{muthetal} 
\een  

We have from \cite{Klartag} Lemma 2.1 and Remarks page 14 the following result.   

\begin{prop}\label{lcritic}  
1) The function  
\be 
r\in (0, +\infty)\to \varphi(r,\theta,\l):=
\frac{\|\A(r\theta+\l)-\y\|_2^2}{2}+\|r\theta+\l\|_1-(p-1)\ln(r)+h(0)
\ee 
is convex and its unique critical point $r(\theta,\l)$
is the mode of (\ref{muthetal}).  

2) By denoting $M(\theta,\l)=\exp\Big(-\varphi(r(\theta,\l),\theta,\l)\Big)$, we obtain 
 \be 
\frac{M(\theta,\l)r(\theta,\l)}{p}\leq J_p(\theta,\l)\leq \frac{M(\theta,\l)r(\theta,\l)(p-1)!\exp(p-1)}{(p-1)^{p}},
\ee 
and for $q >0$, 
\be 
\int_{qr(\theta,\l)}\exp(-\varphi(r,\theta,\l))dr\leq 
\frac{p\Gamma(p,(p-1)q)\exp(p-1)}{(p-1)^p}\int_0^{+\infty}\exp(-\varphi(r,\theta,\l))dr.   
\ee

3) We have 
\be  
\frac{|S|\inf_{\theta\in S}M(\theta,\l)r(\theta,\l)}{p}\leq Z_f\leq \frac{|S|(p-1)!\exp(p-1)}{(p-1)^{p}}\sup_{\theta\in S}M(\theta,\l)r(\theta,\l). 
\ee 

4) if $\x$ is drawn from the distribution $c$ (\ref{c}), 
then $\|\x-\l\|_2\leq qr(\theta,\l)$ with the probability at least equal to 
$1-\frac{p\Gamma(p,(p-1)q)\exp(p-1)}{(p-1)^p}$.  
\end{prop}
\subsection{Calculation of the mode of (\ref{muthetal}) and the partition function 
(\ref{Jp})}  
Now, we are going to calculate the mode $r(\theta,\l)$, and 
the partition function $J_p(\theta,\l)$. 
The calculations are similar to the case LASSO=\{0\}, but  
we need new notations. 
The vector $\y_{\l}=\y-\A\l$, $\s_{\l}=\cos(\theta_{\l})$ where 
$\theta_{\l}$ denotes the angle $(\A\theta,\y_{\l})$, 
$b_{\l}=\|\y_{\l}\|_2s_{\l}$. 
The components of the vector 
$\l\in LASSO$ are denoted by $l_1$, \ldots, $l_p$.  
For $\theta\in\Rb^p$, we set 
\be 
S_0&=&\{i\in\{1, \ldots, p\}:\quad \theta_i=0\},\\
S_+&=&\{i\in\{1, \ldots, p\}:\quad \theta_i\neq 0,\quad  \theta_i l_i \geq 0\},\\
S_-&=&\{i\in\{1, \ldots, p\}:\quad \theta_i\neq 0,\quad  \theta_i l_i < 0\}. 
\ee 
The cardinality of $S_-$ is denoted by $|S_-|$, and the order 
statistic of the sequence 
$\frac{|l_i|}{|\theta_i|}$, for $i\in S_-$ is denoted by  
\be 
l\theta(0):=0 \leq l\theta(1):=\frac{|l|}{|\theta|}(1)\leq \ldots  \leq l\theta(|S_-|):=\frac{|l|}{|\theta|}(|S_-|)\leq 
l\theta(|S_-|+1):=+\infty.  
\ee 
Using these new notations, 
we obtain 
\be 
\|r\theta+\l\|_1=\sum_{i\in S_0}|l_i|+\sum_{i\in S_+}|\theta_i|(r+\frac{l_i}{\theta_i})+
\sum_{i\in S_-}|\theta_i||r-\frac{|l_i|}{|\theta_i|}|. 
\ee 
If $l\theta(k)\leq r < l\theta(k+1)$, then  
\be 
\|r\theta+\l\|_1=\|\theta\|_{1,k}r+c_k, 
\ee 
where 
\be 
c_k:=\sum_{i\in S_0}|l_i|+\sum_{i\in S_+}|l_i|-\sum_{i=0, i\in S_-}^k |l(i)|
+\sum_{i=k+1}^{|S_-|}|l(i)|,\\
\|\theta\|_{1,k}:=
\sum_{i\in S_+}|\theta_i|+\sum_{i=0, i\in S_-}^k|\theta(i)|-\sum_{i=k+1,i\in S_-}^{|S_-|}|\theta(i)|.
\ee 
Observe that $\|\theta\|_{1,|S_-|}=\|\theta\|_1$, and if $\A\theta=0$ then 
\be 
J_p(\theta,\l)=\exp(-\frac{\|\y_{\l}\|_2^2}{2})\sum_{k=0}^{|S_-|}\int_{l\theta(k)}^{l\theta(k+1)}\exp(-\|\theta\|_{1,k}r)r^{p-1}dr.
\ee 
Now, we have the following. 

\begin{prop} If $\A\theta=0$, then 
\ben 
\exp(\frac{\|\y_l\|_2^2}{2})J_p(\theta,\l)=\sum_{k\in I_1(\theta)}\frac{\exp(-c_k)\Big((l\theta(k+1))^p-(l\theta(k))^p\Big)}{p}+\nonumber\\
\sum_{k\in I_2(\theta)}\frac{\exp(-c_k)\Big(\Gamma(p,l\theta(k))-\Gamma(p,l\theta(k+1))\Big)}{\|\theta\|_{1,k}^p},\label{Athetazerol}
\een 
where 
\be 
I_1(\theta)=\{k\in\{0,\ldots, |S_-|\},\,\,\mbox{such that}\,\,\|\theta\|_{1,k}=0\},\\ 
I_2(\theta)=\{k\in\{0,\ldots, |S_-|\},\,\,\mbox{such that}\,\,\|\theta\|_{1,k}> 0\}.
\ee

Now, we are going to give the closed form of $J_p(\theta,\l)$ 
when $\A\theta\neq 0$. 
\end{prop}
We observe for $l\theta(k)\leq r < l\theta(k+1)$ that 
\be 
\frac{\|r\A\theta+\A\l-\y\|_2^2}{2}+\|r\theta+\l\|_1=\alpha_k+\frac{(\|\A\theta\|_2r+\beta_k)^2}{2},
\ee
where 
\be 
&&\beta_k=\frac{\|\theta\|_{1,k}}{\|\A\theta\|_2}-b_{\l},\\
&&\alpha_k=\frac{\|\A\l-\y\|_2^2-\beta_k^2}{2}+c_k.  
\ee 
Moreover if $k\in I_1(\theta)$, then $\beta_k=-b_{\l}$. 
Observe also that $\beta_{|S_-|}$ is bounded below by $-\|\y_{\l}\|_2\sup(s_{\l}:\,\theta\in S)$. 
It follows that 
\be 
J_p(\theta,\l)=\sum_{k\in I_1(\theta)}\exp(-\alpha_k)J_{p,k}(\theta,\l)+
\sum_{k\in I_2(\theta)}\exp(-\alpha_k)J_{p,k}(\theta,\l),
\ee 
where 
\be 
J_{p,k}(\theta,\l)=\int_{l\theta(k)}^{l\theta(k+1)}\exp(-\frac{(\|\A\theta\|_2r+\beta_k)^2}{2})r^{p-1}dr.
\ee 
The calculation of 
\be 
\int_{l\theta(k)}^{l\theta(k+1)}\exp(-\frac{(\|\A\theta\|_2r+\beta_k)^2}{2})r^{p-1}dr
\ee 
is similar to Proposition (\ref{closedfromlassozero}), and 
depends on the sign of  
\be 
x_k=\|\A\theta\|_2l\theta(k)+\beta_k,\\
y_k=\|\A\theta\|_2l\theta(k+1)+\beta_k. 
\ee 
Let $k_0=\max(k: x_k < 0)$ and $k_1=\min(k: y_k >0)$. 
Observe that $k_0+1\geq k_1$, and then $y_k \leq 0$ for $k< k_1\leq k_0+1$. 
It follows for $k< k_1$ that 
\be 
J_{p,k}(\theta,\l)=\frac{1}{\|\A\theta\|_2^p}\sum_{j=0}^{p-1}2^{\frac{j-1}{2}}\binom{p-1}{j}(-\beta_k)^{p-1}
\Big(\gamma(\frac{j+1}{2},\frac{x_k^2}{2})-\gamma(\frac{j+1}{2},\frac{y_k^2}{2})\Big).
\ee 
If $k_1 < k_0+1$, then $x_{k_1} < 0<y_{k_1}$ and  
\be 
J_{p,k_1}(\theta,\l)=&&\frac{1}{\|\A\theta\|_2^p}\sum_{j=0}^{p-1}2^{\frac{j-1}{2}}\binom{p-1}{j}(-\beta_k)^{p-1-j}\gamma(\frac{j+1}{2},\frac{y_k^2}{2})+\frac{1}{\|\A\theta\|_2^p}\\
&&\sum_{j=0}^{p-1}2^{\frac{j-1}{2}}\binom{p-1}{j}(-\beta_k)^{p-1}\gamma(\frac{j+1}{2},\frac{x_k^2}{2}), 
\ee 
and for $k >k_1$,  
\be 
J_{p,k}(\theta,\l)=&&\frac{1}{\|\A\theta\|_2^p}\sum_{j=0}^{p-1}2^{\frac{j-1}{2}}\binom{p-1}{j}(-\beta_k)^{p-1}
\Big(\gamma(\frac{j+1}{2},\frac{y_k^2}{2})-\gamma(\frac{j+1}{2},\frac{x_k^2}{2})\Big). 
\ee  
If $k_1=k_0+1$, then $0\leq x_{k_1} < y_{k_1}$, and for $k\geq k_1$, 
\be 
J_{p,k}(\theta,\l)=\frac{1}{\|\A\theta\|_2^p}\sum_{j=0}^{p-1}2^{\frac{j-1}{2}}\binom{p-1}{j}(-\beta_k)^{p-1}
\Big(\gamma(\frac{j+1}{2},\frac{y_k^2}{2})-\gamma(\frac{j+1}{2},\frac{x_k^2}{2})\Big).
\ee  

Now we can show that 
$J_p(\theta,\l)$ converges to (\ref{Athetazerol}) as $\A\theta\to 0$, 
and we obtain 
an approximation similar to (\ref{gammaexpansion}) for $J_{p,k}(\theta,\l)$ as $\A\theta\to 0$ for 
each $k\in I_1(\theta)$. 


\section{MCMC diagnosis}
Here we take $p=7$, $n=4$, $\A \sim \Bc(\pm \frac{1}{\sqrt{n}})$ and 
for simplicity we consider $\y=0$. 
We sample from the distribution $c$ (\ref{c})
using Hastings-Metropolis algorithm  $(\x^{(t)})$ and propose the test 
$\|\x^{(t)}\|_2\leq q r(\theta^{(t)})$ 
as a criterion for the convergence. Here 
$\theta^{(t)}:=\frac{\x^{(t)}}{\|\x^{(t)}\|_2}$. 
We recall that if $\x$ is drawn from the target distribution $c$, then 
$\|\x\|_2\leq q r(\theta)$ with the probability at least equal to 
$P(q,p)$. Table 2 gives the values of the probability $P(q, p)$. Note that for $q \geq 2.5$ the criterion $\|\x^{(t)}\|_2\leq q r(\theta^{(t)})$ is satisfied with a large probability.

\renewcommand{\arraystretch}{0.9} 
\setlength{\tabcolsep}{0.07cm} 
 \begin{table}[!ht]
\begin{center}
\begin{tabular}{|c|c|c|c|c|c|c|c|}
\hline $q$ &2 &
 2.5 &
 3 &
   3.5 &
   4 &
   4.5&
   5\\
\hline $P(q,p)$ & 0.6672&
0.9446  &
0.9924 &
 0.9991&
 0.9999  &
 1.0000 &
 1.0000  \\
\hline  
\end{tabular}
\caption{ Values of the probability $P(q, p)$ for $p=7$.}
\end{center}
\end{table}

\subsection{Independent sampler (IS)}  
The proposal distribution 
\be 
Q(\x_2,\x_1)=p(\x_2)=\frac{1}{2^p}\exp(-\|\x_2\|_1),\quad \forall\,\x_1, \x_2.
\ee 
The ratio   
\be 
\frac{c(\x)}{p(\x)}\leq \frac{2^p}{Z},\quad \forall\,\x.
\ee 
It's known that MCMC $(x^{(t)})$ with the target distribution $c$ and the proposal distribution $p$ 
is uniformly ergodic \cite{Mengersen}:  
\be 
\sup_{A\subset \mathcal{B}(\Rb^p)}|\Pb(\x^{(t)}\in A\,|\,\x^{(0)})-\int_{A}c(\x)d\x|\leq
(1-\frac{Z}{2^p})^t. 
\ee
Here $Z\approx 2.2142$ and then $(1-\frac{Z}{2^p})=0.9827$. 
Figure 4(a) shows respectively the plot of $t\to 5r(\theta^{(t)})$ and $t\to \|\x^{(t)}\|_2$. 

\subsection{Random-walk (RW) Metropolis algorithm} 
We do not know if the target distribution $c$
satisfies the curvature condition in \cite{Roberts} Section 6. 
Here we propose  
to analyse the convergence of the Random walk Metropolis algorithm 
$(\x^{(t)})$ using the criterion $\|\x^{(t)}\|_2\leq q r(\theta^{(t)})$. 
Figure 4(b) shows respectively the plot of $t\to 5 r(\theta^{(t)})$ and $t\to \|\x^{(t)}\|_2$.  

Figures 4 show that contrary to independent sampler algorithm, 
the random walk (RW) algorithm satisfies early  
the criterion $\|\x^{(t)}\|_2\leq 5 r(\theta)$. More precisely   
\begin{itemize}
\item[1)] the independent sampler (IS) algorithm 
begins to satisfy the criterion $\|\x^{(t)}\|_2\leq 5 r(\theta^{(t)})$
at $t=8\times10^5$ iteration.

\item[2)] The RW algorithm begins to satisfy the criterion 
$\|\x^{(t)}\|_2\leq 3.5 r(\theta^{(t)})$
at $t=939065$ iteration, but the IS algorithm 
never satisfies the criterion $\|\x^{(t)}\|_2\leq 3.5 r(\theta^{(t)})$. 
\end{itemize}
We finally compare IS and RW algorithms using 
the fact that $\int_{\R^p} \x c(\x)d\x=0$. 
The best algorithm will furnish 
the best approximation of the integral $\int_{\R^p} \x c(\x)d\x$.
Table 3 gives the estimators 
$\frac{1}{N} \sum_{t=1}^{N} x_{IS}^{(t)}\approx \int_{\R^p} \x c(\x)d\x$ and 
$\frac{1}{N} \sum_{t=1}^{N} x_{RW}^{(t)}\approx \int_{\R^p} \x c(\x)d\x$. 
It follows that $\|\frac{1}{N} \sum_{t=1}^{N} x_{IS}^{(t)}\|_2 = 0.0187$ and $\|\frac{1}{N} \sum_{t=1}^{N} x_{RW}^{(t)}\|_2 = 0.0041$. We conclude that the random walk algorithm wins 
for both criteria against independent sampler algorithm.   

\renewcommand{\arraystretch}{0.9} 
\setlength{\tabcolsep}{0.07cm} 
 \begin{table}[!ht]
\begin{center}
\begin{tabular}{|c|c|c|c|c|c|c|c|}
\hline  &$x_1$ &
 $ x_2$ &
   $x_3$&
   $x_4$&
   $x_5$&
   $x_6$&
   $x_7$\\
\hline $\hat{x}_{IS}$ &-0.0005&
  -0.0037&
 0.0016&
 0.0164&
   0.0050&
  0.0021&
  -0.0058 \\
\hline  
 $\hat{x}_{RW}$&0.0005 &
-0.0019&
 -0.0002&
 0.0012&
  -0.0005&
 0.0031&
  -0.0011 \\
\hline 
\end{tabular}
\caption{ $N=10^6$, $p=7$, $n=4$ and $q=5$.}
\end{center}
\end{table}

\begin{figure}[!ht]\centering
\includegraphics[width=16cm]{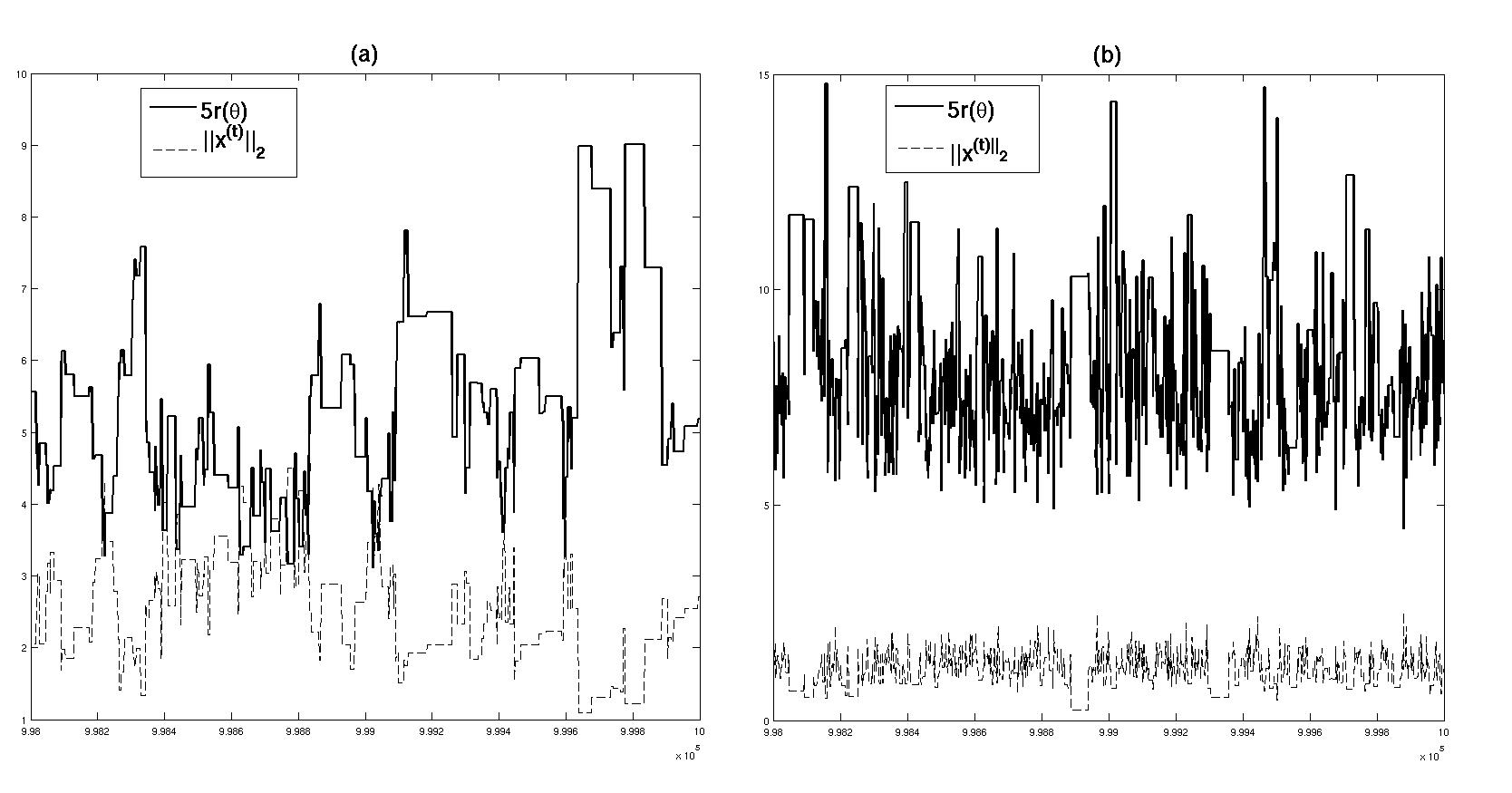} 
\caption{(a): Test of convergence of MCMC algorithm with proposal distribution $p(\x_2)$. (b): Test of convergence of MCMC algorithm with $\Nc(0, 0.5\I_p)$ proposal distribution. $N=10^6$ iterations, $p=7$, $n=4$, $q=5$, $\A \sim \Bc(\pm\frac{1}{\sqrt{n}})$, $\y=0$ and $\l=0$.}
\end{figure}    

\section{Conclusion}
We studied the geometry of bayesian LASSO using polar coordinates  
and calculated the partition function. 
We obtained a concentration inequality and derived MCMC convergence diagnosis
for the convergence of hasting metropolis algorithm. 
We showed that the random walk MCMC with the variance 0.5 wins again 
the independent sampler with the Laplace proposal distribution.

\end{document}